\documentstyle[12pt]{amsart}
\newtheorem{theorem}{Theorem}[section]
\newtheorem{corollary}[theorem]{Corollary}
\newtheorem{remark}[theorem]{Remark}
\newtheorem{lemma}[theorem]{Lemma}

\numberwithin{equation}{section}

\begin{document}

\title[ Long time existence of the symplectic mean curvature flow]
{Long time existence of the symplectic mean curvature flow}

\author{Xiaoli Han, Jiayu Li}

\address{Xiaoli Han, Department of Mathematical Sciences, Tsinghua University \\ Beijing 100084, P. R. of China.}
\email{xlhan@@math.tsinghua.edu.cn}

\address{Jiayu Li, School of Mathematical Sciences, University of Science and Technology of China Hefei 230026 \\ AMSS CAS Beijing 100190, P. R. China}
\email{lijia@@amss.ac.cn}

\thanks {The first author was supported by NSF in China, No.10901088. The second author was supported by NSF in China, No. 11071236.}

\begin{abstract}
Let $(M,\overline{g})$ be a K\"ahler surface with a constant holomorphic sectional curvature $k>0$,
and $\Sigma$ an immersed symplectic surface in $M$. Suppose $\Sigma$ evolves along the mean curvature flow in $M$.
In this paper, we show that the symplectic mean curvature flow exists for long time and converges to a holomorphic curve if
the initial surface satisfies  $|A|^2\leq \frac{2}{3}|H|^2+\frac{1}{2}k$ and $\cos\alpha\geq \frac{\sqrt{30}}{6}$
or $|A|^2\leq \frac{2}{3}|H|^2+\frac{4}{5}k\cos\alpha$ and $\cos\alpha\geq\frac{251}{265}$.

\end{abstract}
\maketitle

\section{Introduction}
Let $M$ be a K\"ahler surface. Let $\omega$ be
the K\"ahker form on $M$ and let $J$ be a complex structure
compatible with $\omega$. The Riemannian metric $\overline{g}$
on $M$ is defined by
$$
\overline{g}( U,V ) =\omega(U,JV).
$$
For a compact oriented real surface $\Sigma$ which is smoothly
immersed in $M$, the k\"ahler angle \cite{CW} $\alpha$ of $\Sigma$ in $M$ was defined by

$$\omega|_\Sigma=\cos\alpha d\mu_\Sigma$$ where $d\mu_\Sigma$ is
the area element of $\Sigma$ in the induced metric from
$\overline{g}$. We say that $\Sigma$ is a symplectic surface if
$\cos\alpha > 0$. The problem is whether one can deform a
symplectic surface to a holomorphic curve ($\cos\alpha \equiv 1$)
in a K\"ahler surface. One way is to use mean curvature flows
(c.f. \cite{CL1}, \cite{W1}, and \cite{HL}), the other way is to
use variational method \cite{HL1}. Chen-Tian \cite{CT}, Chen-Li \cite{CL1} and Wang
\cite{W1} proved that, in a K\"ahler-Einstein surface, if the
initial surface is symplectic, then along the mean curvature flow,
at every time $t$ the surface $\Sigma_t$ is symplectic, which we
call a symplectic mean curvature flow. They also showed that,
there is no type I singularity along a symplectic mean curvature
flow. The symplectic mean curvature flow exists globally and
converges at infinity in graphic cases (c.f. \cite{CLT}, and
\cite{W2}). Han-Li \cite{HL} proved that, in a K\"ahler-Einstein
surface with positive scalar curvature, if the initial surface is
sufficiently close to a holomorphic curve, the symplectic mean
curvature flow exists globally and converges to a holomorphic
curve at infinity. The second type singularity was also studied by
Chen-Li \cite{CL2}, Han-Li \cite{HL2}, Neves \cite{N}, Neves-Tian
\cite{NT}, etc.

Even though one thinks the mean curvature flows may produce
minimal surfaces, there are rather few results on the global
existence and convergence to a minimal surface at infinity of the
mean curvature flows.

In this paper, we consider the case that $M={\bf CP}^2$, i.e.
$(M,\overline{g})$ is a K\"ahler surface with constant
holomorphic sectional curvature $k>0$. We find the condition that
$|A|^2\leq \frac{2}{3}|H|^2+\frac{1}{2}k$ and $\cos\alpha\geq
\frac{\sqrt{30}}{6}$ or $|A|^2\leq \frac{2}{3}|H|^2+\frac{4}{5}k\cos\alpha$
and $\cos\alpha\geq\frac{251}{265}$ is preserved by the mean curvature flow, and
consequently, we show that the symplectic mean curvature flow
exists for long time and converges to a holomorphic curve at
infinity if the initial surface satisfies one of the conditions. As we
know that it is the first long time existence and convergence
result without graphic structure or small initial data conditions.
The main point is to find the pinching condition in our theorem,
which was inspired by Andrews-Baker \cite{AB} and Huisken
\cite{Hu1}.

We believe that, {\it the symplectic mean curvature flow exists
globally and converges to a holomorphic curve at infinity in a
K\"ahler-Einstein surface with positive scalar curvature.}

{\bf Acknowledgement:} The authors thank Jun Sun, Liuqing Yang and Chao Wu for their valuable suggestions.

\section{Preliminaries}

Suppose that $\Sigma$ is submanifold in a Riemannian manifold $M$,
we choose an orthonormal basis $\{e_i\}$ for $T\Sigma$ and $\{e_\alpha\}$ for $N\Sigma$. Recall the evolution equation for the second fundamental form
$h^\alpha_{ij}$ and $|A|^2$ along the mean curvature flow, (see \cite{CL1}, \cite{Sm1}, \cite{W1})
\begin{lemma}
For a mean curvature flow $F: \Sigma \times [0,t_{0}) \to M$, the
second fundamental form $h^{\alpha}_{ij}$ satisfies the following
equation
\begin{eqnarray}\label{e1}
\frac{\partial}{\partial t}h^{\alpha}_{ij}&=&
    \Delta h^{\alpha}_{ij} +(\overline{\nabla}_{\partial_{k}}K)_{\alpha ijk}+(\overline{\nabla}_{\partial_{j}}K)_{\alpha kik} \nonumber\\
    && -2K_{lijk}h^{\alpha}_{lk}+ 2K_{\alpha\beta jk}h^{\beta}_{ik}+ 2K_{\alpha\beta ik}h^{\beta}_{j k} \nonumber\\
    && -K_{lkik}h^{\alpha}_{lj}-K_{lkjk}h^{\alpha}_{il}+K_{\alpha k\beta k}h^{\beta}_{ij} \nonumber\\
    && -H^{\beta}(h^{\beta}_{ik}h^{\alpha}_{jk}+h^{\beta}_{jk}h^{\alpha}_{ik}) \nonumber\\
    && +h^{\alpha}_{im}h^{\beta}_{mk}h^{\beta}_{kj}-2h^{\beta}_{im}h^{\alpha}_{mk}h^{\beta}_{kj}
                    +h^{\beta}_{ik}h^{\beta}_{km}h^{\alpha}_{mj} \nonumber\\
    && +h^{\alpha}_{km}h^{\beta}_{mk}h^{\beta}_{ij}+h^{\beta}_{ij}<e_{\beta},\overline{\nabla}_{H}e_{\alpha}>,
\end{eqnarray}
where $K_{ABCD}$ is the curvature tensor of $M$ and
$\overline{\nabla}$ is the covariant derivative of $M$. Therefore
\begin{eqnarray}\label{e2}
\frac{\partial}{\partial t}|A|^{2}&=&
     \Delta |A|^{2} -2|\nabla A|^{2}
       +[(\overline{\nabla}_{\partial_{k}}K)_{\alpha ijk}+(\overline{\nabla}_{\partial_{j}}K)_{\alpha kik}]h^{\alpha}_{ij} \nonumber\\
    && -4K_{lijk}h^{\alpha}_{lk}h^{\alpha}_{ij} +8K_{\alpha\beta jk}h^{\beta}_{ik}h^{\alpha}_{ij}
       -4K_{lkik}h^{\alpha}_{lj}h^{\alpha}_{ij}+ 2K_{\alpha k\beta k}h^{\beta}_{ij}h^{\alpha}_{ij} \nonumber\\
    && +2\sum_{\alpha,\beta,i,j}(\sum_{k}(h^{\alpha}_{ik}h^{\beta}_{jk}-h^{\alpha}_{jk}h^{\beta}_{ik}))^{2}
       +2\sum_{\alpha,\beta}(\sum_{ij}h^{\alpha}_{ij}h^{\beta}_{ij})^{2}.
\end{eqnarray}
\end{lemma}

\begin{corollary}
Along the mean curvature flow, the mean curvature vector satisfies
\begin{eqnarray}\label{e3}
\frac{\partial}{\partial t}|H|^{2}&=&
    \Delta |H|^{2} -2|\nabla H|^{2}+2K_{\alpha k\beta k}H^{\alpha}H^{\beta}
               +2\sum_{ij}(\sum_{\alpha} H^{\alpha}h^{\alpha}_{ij})^2.
\end{eqnarray}

\end{corollary}

Suppose that $M$ is a compact K\"ahler surface.
Let $\Sigma$ be a smooth surface in $M$.
The K\"ahler angle of $\Sigma$ in $M$ is defined by (\cite{CW})
$$ \omega|_\Sigma=\cos\alpha d\mu_\Sigma
$$ where $d\mu_\Sigma$ is the area element of $\Sigma$ of the
induced metric from $\overline{g}$. We call $\Sigma$ a {\it
symplectic} surface if $\cos\alpha>0$, a {\it Lagrangian} surface if
$\cos\alpha\equiv 0$, a {\it holomorphic curve} if $\cos\alpha\equiv
1$. Recall the evolution equation of $\cos\alpha$ (\cite{CL1}, \cite{W1}),
\begin{lemma} Along the mean curvature flow, $\cos\alpha$ satisfies
\begin{eqnarray}\label{e6}
(\frac{\partial}{\partial t}-\Delta)\cos\alpha=|\overline{\nabla}J_{\Sigma_t}|^2\cos\alpha+Ric(Je_1, e_2)\sin^2\alpha.
\end{eqnarray} where $|\overline{\nabla}J_{\Sigma_t}|^2=|h^3_{1k}-h^4_{2k}|^2+|h^3_{2k}+h^4_{1k}|^2$, $\{e_1, e_2, e_3, e_4\}$ is
any orthonormal basis for $TM$ such that $\{e_1, e_2\}$ is the basis for $T\Sigma$ and $\{e_3, e_4\}$ is the basis for $N\Sigma$.
\end{lemma}
It is proved in \cite{CL1} and \cite{HL1} that
 \begin{eqnarray}\label{e17}
|\overline{\nabla}J_{\Sigma_t}|^2\geq \frac{1}{2}|H|^2\end{eqnarray} and
\begin{eqnarray}\label{e18}|\nabla\cos\alpha|^2\leq\sin^2\alpha|\overline{\nabla}J_{\Sigma_t}|^2. \end{eqnarray}

Now suppose $M$ is a K\"ahler surface with constant holomorphic sectional curvature $k$, then from Theorem $2.1$ and Theorem $2.3$ in \cite{Y},
we have
\begin{lemma}
 $M$ has a curvature tensor of the form
\begin{eqnarray}\label{e4}
K_{kjih}=-\frac{k}{4}[(g_{kh}g_{ji}-g_{jh}g_{ki})+(J_{kh}J_{ji}-J_{jh}J_{ki})-2J_{kj}J_{ih}].
\end{eqnarray} Thus $M$ is symmetric. Furthermore,
$M$ is Einstein
\begin{eqnarray}\label{e5}
K_{ji}=\frac{3}{2}k \overline{g}_{ij}.
\end{eqnarray}
\end{lemma}

\section{Pinching estimate}

In this section we want to show our pinching inequality is preserved by the symplectic mean curvature flow. Before proving our theorem,
we deduce the local expression of the complex structure of the K\"ahler surface.
Let $M$ be a K\"ahler surface with the K\"ahler metric $\bar{g}$ and $\Sigma$ be a real surface in $M$. Suppose  $\omega$
is the associated K\"ahler form and $J$ is the complex structure compatible with $\bar{g}$ and $\omega$, i.e,
$$\omega(X, Y)=\bar{g}(JX, Y)=\langle JX, Y\rangle$$ for any $X, Y\in TM$. Fix $p\in M$. We choose the local frame of
$M$ around $p$ $\{e_1, e_2, e_3, e_4\}$ such that $\{e_1, e_2\}$ is the frame of the tangent bundle $T\Sigma$ and $\{e_3, e_4\}$ is the frame
of the normal bundle $N\Sigma$. Suppose
$$Je_1=xe_2+ye_3+ze_4.$$ Then using $\langle JX, Y\rangle=-\langle X, JY\rangle$ and
$$-e_1=xJe_2+yJe_3+zJe_4,$$ we have
\begin{eqnarray*}
\left \{ \begin{array}{clcr} x\langle Je_2, e_4\rangle+y\langle Je_3, e_4\rangle= &0\\
x\langle Je_2, e_3\rangle-z\langle Je_3, e_4\rangle= &0\\
-y\langle Je_2, e_3\rangle-z\langle Je_2, e_4\rangle= &0 \end{array}\right..
\end{eqnarray*}
Suppose $y\neq 0$. Set $\langle Je_2, e_4\rangle=A$.  Then we have
$$\langle Je_2, e_3\rangle=-\frac{z}{y}A,~~~~~~~~~~~~~~~~~~~~~~~~~~~~\langle Je_3, e_4\rangle=-\frac{x}{y}A.$$
Thus
\begin{eqnarray*}
J=\left (\begin{array}{clcr} 0 &x &y &z \\
-x &0 &-\frac{z}{y}A &A\\
-y &\frac{z}{y}A &0 &-\frac{x}{y}A\\
-z &-A &\frac{x}{y}A &0 \end{array}\right).
\end{eqnarray*} As $J$ is isometric, we have
\begin{eqnarray*}
\left \{\begin{array}{clcr}  x^2+y^2+z^2 &=1 \\
x^2+(\frac{z}{y})^2A^2+A^2 &=1 \end{array}\right.,
\end{eqnarray*} we can obtain that $A^2=y^2$, i.e., $A=\pm y$. Thus we see that
\begin{eqnarray}\label{e14}
J=\left (\begin{array}{clcr} 0 &x &y &z \\
-x &0 &-z &y\\
-y &z &0 &-x\\
-z &-y &x &0 \end{array}\right),
\end{eqnarray} or
\begin{eqnarray}\label{e15}
J=\left (\begin{array}{clcr} 0 &x &y &z \\
-x &0 &z &-y\\
-y &-z &0 &x\\
-z &y &-x &0 \end{array}\right).
\end{eqnarray}

If $y=0$, then by the same argument we can see that $J$ also has the form (\ref{e14}) or (\ref{e15}).
By the definition of the K\"ahler angle, we know that
$$x=\cos\alpha=\omega(e_1, e_2)=\langle Je_1, e_2\rangle.$$ If we assume the K\"ahler form is anti-self-dual,
then $J$ has the form (\ref{e15}).

We begin by estimating the gradient terms:

\begin{lemma}\label{L2}
For any $\eta>0$ we have the inequality
\begin{equation}\label{e27}
|\nabla A|^2\geq(\frac{3}{n+2}-\eta)|\nabla H|^2-\frac{2}{n+2}(\frac{2}{n+2}\eta^{-1}-\frac{n}{n-1})|w|^2,
\end{equation} where $w_i^\alpha=\sum_{l} K_{\alpha lil}$, $|w^\alpha|^2=\sum_{i}|w^\alpha_i|^2$ and $|w|^2=\sum_{\alpha}|w^\alpha|^2$ .
\end{lemma}

{\it Proof.} Similar as \cite{Ha} and \cite{Hu1} we decompose the tensor $\nabla A$ into
$$\nabla_i h^\alpha_{jk}=E^\alpha_{ijk}+F^\alpha_{ijk},$$
where
\begin{eqnarray*}
E^\alpha_{ijk} &=&\frac{1}{n+2}(\nabla_i H^\alpha\cdot g_{jk}+\nabla_j H^\alpha\cdot g_ik+\nabla_k H^\alpha\cdot g_{ij})\\
&&-\frac{2}{(n+2)(n-1)}w_i^\alpha g_{jk}+\frac{n}{(n+2)(n-1)}(w_j^\alpha g_{ik}+w_k^\alpha g_{ij}).
\end{eqnarray*}
It is easy to get that, $\langle E^\alpha_{ijk}, F^\alpha_{ijk}\rangle=0.$ Furthermore,
\begin{eqnarray*}
|E^\alpha|^2 &=&\frac{3}{n+2}|\nabla H|^2+\frac{2n}{(n+2)(n-1)}|w^\alpha|^2+\frac{4}{n+2}\langle w_i^\alpha, \nabla_i H^\alpha\rangle\\
&\geq&(\frac{3}{n+2}-\eta)|\nabla H^\alpha|^2-\frac{2}{n+2}(\frac{2}{n+2}\eta^{-1}-\frac{n}{n-1})|w^\alpha|^2.
\end{eqnarray*} We finish the proof of the Lemma.
\hfill Q. E. D.

\begin{theorem}\label{L1}
Suppose $M$ is a K\"ahler surface with constant holomorphic sectional curvature $k>0$ and $\Sigma$ is a symplectic surface in $M$. Assume that
$|A|^2\leq \frac{2}{3}|H|^2+\frac{1}{2}k$ and $\cos\alpha\geq \frac{\sqrt{30}}{6}$ holds on the initial surface, then it remains true along the symplectic
mean curvature flow.
\end{theorem}

{\it Proof.} From  (\ref{e6}) and (\ref{e5}), we know that
\begin{eqnarray*}
(\frac{\partial}{\partial t}-\Delta)\cos\alpha=|\overline{\nabla}J_{\Sigma_t}|^2\cos\alpha+\frac{3k}{2}\cos\alpha\sin^2\alpha.
\end{eqnarray*} Thus at any time $t$, $\cos\alpha\geq \frac{\sqrt{30}}{6}$ if it holds on the initial surface.

Since $M$ is symmetric, by (\ref{e2}) we know that
\begin{eqnarray}\label{e7}
\frac{\partial}{\partial t}|A|^{2}&=&
     \Delta |A|^{2} -2|\nabla A|^{2}\nonumber\\
    && -4K_{lijk}h^{\alpha}_{lk}h^{\alpha}_{ij} +8K_{\alpha\beta jk}h^{\beta}_{ik}h^{\alpha}_{ij}
       -4K_{lkik}h^{\alpha}_{lj}h^{\alpha}_{ij}+ 2K_{\alpha k\beta k}h^{\beta}_{ij}h^{\alpha}_{ij} \nonumber\\
    && +2\sum_{\alpha,\beta,i,j}(\sum_{k}(h^{\alpha}_{ik}h^{\beta}_{jk}-h^{\alpha}_{jk}h^{\beta}_{ik}))^{2}
       +2\sum_{\alpha,\beta}(\sum_{ij}h^{\alpha}_{ij}h^{\beta}_{ij})^{2}.\nonumber
\end{eqnarray}
Now in our case the first four terms reduce to
\begin{eqnarray*}
-4K_{lijk}h^{\alpha}_{lk}h^{\alpha}_{ij}&=&-4K_{1212}(h^\alpha_{12})^2-4K_{1221}h^\alpha_{11}h^\alpha_{22}\\
&&-4K_{2112}h^\alpha_{11}h^\alpha_{22}-4K_{2121}(h^\alpha_{12})^2\\&=&-4K_{1212}(2(h^\alpha_{12})^2-2h^\alpha_{11}h^\alpha_{22})\\
&=&-4K_{1212}(|A|^2-|H|^2),
\end{eqnarray*}
and
\begin{eqnarray*}
8K_{\alpha\beta jk}h^{\beta}_{ik}h^{\alpha}_{ij}&=&8K_{3412}h^3_{i1}h^4_{i2}+8K_{3421}h^3_{i2}h^4_{i1}\\&&+
8K_{4312}h^4_{i1}h^3_{i2}+8K_{4321}h^4_{i2}h^3_{i1}\\&=&16K_{1234}(h^3_{1i}h^4_{2i}-h^3_{i2}h^4_{1i})\\&=& 8K_{1234}(|A|^2-|\overline{\nabla} J_{\Sigma_t}|^2),
\end{eqnarray*} and
\begin{eqnarray*}
-4K_{lkik}h^{\alpha}_{lj}h^{\alpha}_{ij}&=&-4K_{1212}(h^\alpha_{1j})^2-4K_{2121}(h^\alpha_{2j})^2\\ &=& -4K_{1212}|A|^2,
\end{eqnarray*} and
\begin{eqnarray*}
2K_{\alpha k\beta k}h^{\beta}_{ij}h^{\alpha}_{ij}&=&2K_{3k3k}(h^3_{ij})^2+2K_{4k4k}(h^4_{ij})^2+4K_{3k4k}h^3_{ij}h^4_{ij}\\ &=&
2K_{33}(h^3_{ij})^2-2K_{3434}(h^3_{ij})^2\\&&+2K_{44}(h^4_{ij})^2-2K_{3434}(h^4_{ij})^2+4K_{34}h^3_{ij}h^4_{ij}\\
&=&3k|A|^2-2K_{3434}|A|^2,
\end{eqnarray*} where we have used the equality (\ref{e5}). Therefore,
\begin{eqnarray}\label{e8}
\frac{\partial}{\partial t}|A|^{2}&=&
     \Delta |A|^{2} -2|\nabla A|^{2}\nonumber\\
    && +8(K_{1234}-K_{1212})|A|^2+3k|A|^2-2K_{3434}|A|^2 \nonumber\\
&&+4K_{1212}|H|^2-8K_{1234}|\overline{\nabla} J_{\Sigma_t}|^2\nonumber\\
    && +2\sum_{\alpha,\beta,i,j}(\sum_{k}(h^{\alpha}_{ik}h^{\beta}_{jk}-h^{\alpha}_{jk}h^{\beta}_{ik}))^{2}
       +2\sum_{\alpha,\beta}(\sum_{ij}h^{\alpha}_{ij}h^{\beta}_{ij})^{2}.
\end{eqnarray} Similarly, the evolution equation of $|H|^2$ becomes
\begin{eqnarray}\label{e10}
\frac{\partial}{\partial t}|H|^{2}&=&
    \Delta |H|^{2}-2|\nabla H|^{2}+3k|H|^2-2K_{3434}|H|^2\nonumber\\&&
               +2\sum_{ij}(\sum_{\alpha} H^{\alpha}h^{\alpha}_{ij})^2.
\end{eqnarray}

Using  (\ref{e15}) we get that,
\begin{eqnarray}\label{e16}
K_{1212}&=&K_{3434}=\frac{k}{4}(3\cos^2\alpha+1);\nonumber\\
K_{1234} &=&-\frac{k}{4}(z^2+y^2-2x^2)=\frac{k}{4}(3\cos^2\alpha-1).
\end{eqnarray}

Putting (\ref{e16}) into (\ref{e8}) , we get that
\begin{eqnarray*}
\frac{\partial}{\partial t}|A|^{2}&=&
     \Delta |A|^{2} -2|\nabla A|^{2}-k|A|^2-\frac{k}{2}(3\cos^2\alpha+1)|A|^2\\
&&+k(3\cos^2\alpha+1)|H|^2-2k(3\cos^2\alpha-1)|\overline{\nabla} J_{\Sigma_t}|^2 \\ &&
+2\sum_{\alpha,\beta,i,j}(\sum_{k}(h^{\alpha}_{ik}h^{\beta}_{jk}-h^{\alpha}_{jk}h^{\beta}_{ik}))^{2}
       +2\sum_{\alpha,\beta}(\sum_{ij}h^{\alpha}_{ij}h^{\beta}_{ij})^{2}.
\end{eqnarray*} Using the inequality  (\ref{e17}) and
$\cos\alpha\geq \frac{\sqrt{30}}{6}>\frac{\sqrt{3}}{3}$, we obtain that
\begin{eqnarray}\label{e33}
\frac{\partial}{\partial t}|A|^{2}&\leq&
     \Delta |A|^{2} -2|\nabla A|^{2}-k|A|^2-\frac{k}{2}(3\cos^2\alpha+1)|A|^2+2k|H|^2\nonumber\\ &&
+2\sum_{\alpha,\beta,i,j}(\sum_{k}(h^{\alpha}_{ik}h^{\beta}_{jk}-h^{\alpha}_{jk}h^{\beta}_{ik}))^{2}
       +2\sum_{\alpha,\beta}(\sum_{ij}h^{\alpha}_{ij}h^{\beta}_{ij})^{2}.
\end{eqnarray}
Similarly,
\begin{eqnarray}\label{e29}
\frac{\partial}{\partial t}|H|^{2}&=&
    \Delta |H|^{2}-2|\nabla H|^{2}+3k|H|^2-\frac{k}{2}(3\cos^2\alpha+1)|H|^2\nonumber\\&&
               +2\sum_{ij}(\sum_{\alpha} H^{\alpha}h^{\alpha}_{ij})^2.
\end{eqnarray}

Set $Q=|A|^2-\frac{2}{3}|H|^2-bk$. Therefore,\allowdisplaybreaks
\begin{eqnarray}\label{e26}
\frac{\partial}{\partial t}Q &\leq&\Delta Q-2(|\nabla A|^{2}-\frac{2}{3}|\nabla H|^{2})-k|A|^2\nonumber\\ &&
-\frac{k}{2}(3\cos^2\alpha+1)(|A|^2-\frac{2}{3}|H|^2)\nonumber\\&&
+2\sum_{\alpha,\beta,i,j}(\sum_{k}(h^{\alpha}_{ik}h^{\beta}_{jk}-h^{\alpha}_{jk}h^{\beta}_{ik}))^{2}
       +2\sum_{\alpha,\beta}(\sum_{ij}h^{\alpha}_{ij}h^{\beta}_{ij})^{2}\nonumber\\&&-\frac{4}{3}
\sum_{ij}(\sum_{\alpha} H^{\alpha}h^{\alpha}_{ij})^2\nonumber\\ &\leq&  \Delta Q-2(|\nabla A|^{2}-\frac{2}{3}|\nabla H|^{2})-\frac{k}{2}(3\cos^2\alpha+1)Q
\nonumber\\ &&-\frac{bk^2}{2}(3\cos^2\alpha+1)-k|A|^2
\nonumber\\&&
+2\sum_{\alpha,\beta,i,j}(\sum_{k}(h^{\alpha}_{ik}h^{\beta}_{jk}-h^{\alpha}_{jk}h^{\beta}_{ik}))^{2}
       +2\sum_{\alpha,\beta}(\sum_{ij}h^{\alpha}_{ij}h^{\beta}_{ij})^{2}\nonumber\\&&-\frac{4}{3}
\sum_{ij}(\sum_{\alpha} H^{\alpha}h^{\alpha}_{ij})^2.
\end{eqnarray}

First we estimate the gradient terms in (\ref{e26}). In (\ref{e27}) we choose $\eta=\frac{1}{12}, n=2$, then
\begin{eqnarray*}
|\nabla A|^2\geq\frac{2}{3} |\nabla H|^2-2|w|^2,
\end{eqnarray*} where
$$|w|^2=K_{3212}^2+K_{3121}^2+K_{4121}^2+K_{4212}^2.$$
Using (\ref{e15}) again, we obtain that
$$|w|^2=\frac{9k^2}{8}(x^2z^2+x^2y^2)=\frac{9k^2}{8}x^2(1-x^2)=\frac{9k^2}{8}\cos^2\alpha\sin^2\alpha.$$
Thus
\begin{equation}\label{e28}
|\nabla A|^2\geq\frac{2}{3} |\nabla H|^2-\frac{9k^2}{4}\cos^2\alpha\sin^2\alpha.
\end{equation}
In order to estimate the other terms in (\ref{e26}) we do the same way as in \cite{AB}.
Set $R_1=\sum_{\alpha,\beta,i,j}(\sum_{k}(h^{\alpha}_{ik}h^{\beta}_{jk}-h^{\alpha}_{jk}h^{\beta}_{ik}))^{2}$, $R_2=\sum_{\alpha,\beta}(\sum_{ij}h^{\alpha}_{ij}h^{\beta}_{ij})^{2}$, $R_3=
\sum_{ij}(\sum_{\alpha} H^{\alpha}h^{\alpha}_{ij})^2.$
At the point $|H|\neq 0$,
we choose $\{e_3, e_4\}$ for $N\Sigma$
such that $e_3=H/|H|$ and choose the $\{e_1, e_2\}$ for $T\Sigma$ such that $h^3_{ij}=\lambda_i\delta_{ij}$ .
Set $h^\alpha_{ij}=\mathring{h}^\alpha_{ij}+\frac{1}{2}H^\alpha g_{ij}$, then $\mathring{h}^4_{ij}=h^4_{ij},
\mathring{h}^3_{ij}=h^3_{ij}-\frac{1}{2}|H|g_{ij}$. Since $(h^3_{ij})$ is diagonal, we see $(\mathring{h}^3_{ij})$ is also diagonal.
Set  $(\mathring{h}^3_{ij})=\mathring{\lambda}_i\delta_{ij}$. Denote the norm of $(h^\alpha_{ij}), (\mathring{h}^\alpha_{ij})$
by $|h_\alpha|, |\mathring{h}_\alpha|$ respectively. $R_1, R_2, R_3$ reduce to

\begin{eqnarray*}
R_2=\sum_{\alpha,\beta}(\sum_{ij}h^{\alpha}_{ij}h^{\beta}_{ij})^{2} &=&|\mathring{h}_3|^4+|\mathring{h}_4|^4+|\mathring{h}_3|^2|H|^2+\frac{1}{4}|H|^4
\\&&+2(\sum_{ij}\mathring{h}^3_{ij}\mathring{h}^4_{ij})^2;
\end{eqnarray*}
\begin{eqnarray*}
R_1=\sum_{\alpha,\beta,i,j}(\sum_{k}(h^{\alpha}_{ik}h^{\beta}_{jk}-h^{\alpha}_{jk}h^{\beta}_{ik}))^{2} &=&
2\sum_{i,j}(\sum_{k}(h^{3}_{ik}\mathring{h}^{4}_{jk}-h^{3}_{jk}\mathring{h}^{4}_{ik}))^{2};
\end{eqnarray*}
\begin{eqnarray*}
R_3=\sum_{ij}(\sum_{\alpha} H^{\alpha}h^{\alpha}_{ij})^2 &=&|\mathring{h}_3|^2|H|^2+\frac{1}{2}|H|^4.
\end{eqnarray*} Using the fact that $(\mathring{h}^3_{ij})$ is diagonal , then
we have
\begin{eqnarray*}
(\sum_{ij}\mathring{h}^3_{ij}\mathring{h}^4_{ij})^2 &=&
(\sum_{i}\mathring{\lambda}_{i}\mathring{h}^4_{ii})^2 \\ &\leq& (\sum_{i}\mathring{\lambda_i}^2)(\sum_{i}(\mathring{h}^4_{ii})^2)=|\mathring{h}_3|^2 \sum_i(\mathring{h}^4_{ii})^2,
\end{eqnarray*}
\begin{eqnarray*}
\sum_{i,j}(\sum_{k}(h^{3}_{ik}\mathring{h}^{4}_{jk}-h^{3}_{jk}\mathring{h}^{4}_{ik}))^{2} &=& \sum_{i\neq j}(\lambda_i-\lambda_j)^2(\mathring{h}^4_{ij})^2
\\&=&\sum_{i\neq j}(\mathring{\lambda}_i-\mathring{\lambda}_j)^2(\mathring{h}^4_{ij})^2\\&\leq& \sum_{i\neq j}2(\mathring{\lambda}^2_i+\mathring{\lambda}^2_j)^2(\mathring{h}^4_{ij})^2\\ &\leq&2|\mathring{h}_3|^2\sum_{i\neq j}(\mathring{h}_{ij}^4)^2\\
&=&2|\mathring{h}_3|^2 (|\mathring{h}_4|^2-\sum_{i}(\mathring{h}_{ii}^4)^2),
\end{eqnarray*}
so
\begin{eqnarray*}
(\sum_{ij}\mathring{h}^3_{ij}\mathring{h}^4_{ij})^2+\sum_{i,j}(\sum_{k}(h^{3}_{ik}\mathring{h}^{4}_{jk}-h^{3}_{jk}\mathring{h}^{4}_{ik}))^{2} \leq
2|\mathring{h}_3|^2 |\mathring{h}_4|^2.
\end{eqnarray*}

Therefore,
\begin{eqnarray}\label{e34}
2R_1+2R_2-\frac{4}{3}R_3&\leq& 2|\mathring{h}_3|^4+2|\mathring{h}_4|^4+\frac{2}{3}|\mathring{h}_3|^2|H|^2\nonumber\\ &&-
\frac{1}{6}|H|^4+8|\mathring{h}_3|^2|\mathring{h}_4|^2.
\end{eqnarray}

Using these inequalities together with (\ref{e28}), we obtain that
\begin{eqnarray*}
\frac{\partial}{\partial t}Q &\leq&\Delta Q+\frac{9k^2}{2}\cos^2\alpha\sin^2\alpha
-\frac{k}{2}(3\cos^2\alpha+1)Q\\&&-\frac{bk^2}{2}(3\cos^2\alpha+1)-k|A|^2\\ &&+
2 |\mathring{h}_3|^4+2|\mathring{h}_4|^4+\frac{2}{3}|\mathring{h}_3|^2|H|^2-\frac{1}{6}|H|^4+8|\mathring{h}_3|^2|\mathring{h}_4|^2.
\end{eqnarray*} Sine $|H|^2=6(|\mathring{h}_3|^2+|\mathring{h}_4|^2-Q-bk)$, putting it into the above inequality we obtain that
\begin{eqnarray}\label{e12}
\frac{\partial}{\partial t}Q &\leq&\Delta Q-6Q^2\nonumber
\\&&+[8|\mathring{h}_3|^2+12|\mathring{h}_4|^2-12bk+3k-\frac{3k}{2}(\cos^2\alpha+1)]Q\nonumber\\&&
-\frac{bk^2}{2}(3\cos^2\alpha+1)+\frac{9k^2}{2}\cos^2\alpha\sin^2\alpha+3bk^2-6b^2k^2\nonumber\\&&-4|\mathring{h}_4|^4+4k(2b-1)|\mathring{h}_3|^2
+4k(3b-1)|\mathring{h}_4|^2\nonumber\\ &\leq&
\Delta Q-6Q^2\nonumber
\\&&+[8|\mathring{h}_3|^2+12|\mathring{h}_4|^2-12bk+3k-\frac{3k}{2}(\cos^2\alpha+1)]Q\nonumber\\&&
+4k(2b-1)|\mathring{h}_3|^2-(2|\mathring{h}_4|^2-k(3b-1))^2+k^2(3b-1)^2\nonumber\\
&&-\frac{bk^2}{2}(3\cos^2\alpha+1)+\frac{9k^2}{2}\cos^2\alpha\sin^2\alpha+3bk^2-6b^2k^2\nonumber\\ &\leq&
\Delta Q-6Q^2+[8|\mathring{h}_3|^2+12|\mathring{h}_4|^2-12bk+3k-\frac{3k}{2}(\cos^2\alpha+1)]Q\nonumber\\
&&-(2|\mathring{h}_4|^2-k(3b-1))^2+4k(2b-1)|\mathring{h}_3|^2\nonumber\\&&
+k^2(3b^2-\frac{7}{2}b+1)+k^2\cos^2\alpha(\frac{9}{2}\sin^2\alpha-\frac{3}{2}b)\nonumber.
\end{eqnarray}
If $2b-1\leq 0$ and $3b^2-\frac{7}{2}b+1\leq 0$, then $b$ must be equal to $\frac{1}{2}$.
If $\frac{9}{2}\sin^2\alpha-\frac{3}{2}b\leq 0$, then $\sin^2\alpha\leq \frac{1}{6}$, i.e, $\cos^2\alpha\geq
\frac{5}{6}$.

At the point $|H|=0$, we use the following inequality (see \cite{CCK}, \cite{LL}),
\begin{eqnarray}\label{e35}
2\sum_{\alpha,\beta,i,j}(\sum_{k}(h^{\alpha}_{ik}h^{\beta}_{jk}-h^{\alpha}_{jk}h^{\beta}_{ik}))^{2}
       +2\sum_{\alpha,\beta}(\sum_{ij}h^{\alpha}_{ij}h^{\beta}_{ij})^{2} &\leq& 3|A|^4.
\end{eqnarray}
Thus, using (\ref{e28}) and (\ref{e26}) we obtain
\begin{eqnarray*}
\frac{\partial}{\partial t}Q &\leq&\Delta Q+\frac{9k^2}{2}\cos^2\alpha\sin^2\alpha-k|A|^2\\ &&
-\frac{k}{2}(3\cos^2\alpha+1)|A|^2+3|A|^4.\\ &&
\end{eqnarray*} Since $|H|=0$, we have $|A|^2=Q+bk$. Thus,
\begin{eqnarray}\label{e13}
\frac{\partial}{\partial t}Q &\leq&\Delta Q+\frac{9k^2}{2}\cos^2\alpha\sin^2\alpha\nonumber\\ &&
-\frac{3k}{2}(\cos^2\alpha+1)(Q+bk)+3(Q+bk)^2.\nonumber\\ &\leq&
\Delta Q+\frac{9k^2}{2}\cos^2\alpha\sin^2\alpha\nonumber \\&& +[3(|A|^2+bk)-\frac{3k}{2}(\cos^2\alpha+1)]Q\nonumber\\ &&-
\frac{bk^2}{2}(\cos^2\alpha+1)+3b^2k^2\nonumber\\ &\leq& \Delta Q+[3(|A|^2+bk)-\frac{3k}{2}(\cos^2\alpha+1)]Q\nonumber \\
&&+3bk^2(b-\frac{1}{2})+k^2\cos^2\alpha(\frac{9}{2}\sin^2\alpha-\frac{3b}{2}).
\end{eqnarray} Thus we need choose $b\leq\frac{1}{2}$ and $\sin^2\alpha\leq b/3$.

Therefore, we choose $b=\frac{1}{2}$ and $\cos^2\alpha\geq\frac{5}{6}$. Then we have
$$\frac{\partial}{\partial t}Q \leq\Delta Q+CQ.$$ Applying the maximum principle for parabolic equation, we see that
$$Q\leq 0$$ if it holds on initial surface.
\hfill Q. E. D.

\begin{remark}
During reading our paper, Yang Liuqing found the condition that $|A|^2\leq \lambda |H|^2+\frac{2\lambda -1}{\lambda}k$ and $\cos\alpha\geq \sqrt{\frac{7\lambda -3}{3\lambda}}$ ($1/2\leq \lambda\leq 2/3$) is preserved by the symplectic mean curvature flow.
\end{remark}

\section{Long time existence and convergence}
In this section we prove the long time existence of the symplectic mean curvature flow under the assumption of Theorem \ref{L1}.
\begin{theorem}\label{T1}
Under the assumption of Theorem \ref{L1}, the symplectic mean curvature flow exists for long time.
\end{theorem}
{\it Proof.}
Suppose $f$ is a positive function which will be determined later. Now we compute the evolution equation of $\frac{|H|^2}{f(\cos\alpha)}$.
\begin{eqnarray*}
(\frac{\partial}{\partial t}-\Delta)\frac{|H|^2}{f(\cos\alpha)} &=&\frac{(\frac{\partial}{\partial t}-\Delta)|H|^2}{f(\cos\alpha)}
-\frac{|H|^2f'(\frac{\partial}{\partial t}-\Delta)\cos\alpha}{f^2(\cos\alpha)}\\&&+\frac{|H|^2f''|\nabla\cos\alpha|^2}{f^2(\cos\alpha)}
+2\frac{f'\nabla\cos\alpha}{f}\cdot\nabla\frac{|H|^2}{f(\cos\alpha)}.
\end{eqnarray*}
It follows that,
\begin{eqnarray*}
(\frac{\partial}{\partial t}-\Delta)\cos\alpha &=& |\overline{\nabla} J_{\Sigma_t}|^2\cos\alpha+\frac{3}{2}k\sin^2\alpha\cos\alpha\\
&\geq& |\overline{\nabla} J_{\Sigma_t}|^2\cos\alpha.
\end{eqnarray*} By (\ref{e29}) and $\cos^2\alpha\geq\frac{5}{6}$, we have
\begin{eqnarray*}
(\frac{\partial}{\partial t}-\Delta) |H|^2\leq \frac{5k}{4}|H|^2+2|H|^2|A|^2.
\end{eqnarray*}
Putting the above inequality into the evolution equation of $\frac{|H|^2}{f(\cos\alpha)}$, we get that
\begin{eqnarray*}
&&(\frac{\partial}{\partial t}-\Delta)\frac{|H|^2}{f(\cos\alpha)} \leq 2\frac{f'\nabla\cos\alpha}{f}\cdot\nabla\frac{|H|^2}{f(\cos\alpha)}
\\&&+\frac{f(\frac{5k}{4}|H|^2+2|H|^2|A|^2)-|H|^2f'|\overline{\nabla} J_{\Sigma_t}|^2\cos\alpha+|H|^2f''|\nabla\cos\alpha|^2}{f^2(\cos\alpha)}
\\&\leq& 2\frac{f'\nabla\cos\alpha}{f}\cdot\nabla\frac{|H|^2}{f(\cos\alpha)}\\
&&+\frac{f(\frac{5k}{4}|H|^2+2|H|^2(\frac{2}{3}|H|^2+\frac{k}{2})-|H|^2f'|\overline{\nabla} J_{\Sigma_t}|^2\cos\alpha+|H|^2f''|\nabla\cos\alpha|^2}{f^2(\cos\alpha)}\\
&\leq& 2\frac{f'\nabla\cos\alpha}{f}\cdot\nabla\frac{|H|^2}{f(\cos\alpha)}+\frac{9}{4}k\frac{|H|^2}{f(\cos\alpha)}\\&&+
\frac{\frac{4}{3}f|H|^4-\frac{\sqrt{30}}{6}|H|^2f'|\overline{\nabla} J_{\Sigma_t}|^2+\frac{1}{6}|H|^2f''|\overline{\nabla} J_{\Sigma_t}|^2}{f^2(\cos\alpha)},
\end{eqnarray*} where we have used (\ref{e18}), $\sin^2\alpha\leq\frac{1}{6}$ and we assume $f'>0, f''>0$.
We want to find $f$ such that
$$\frac{4}{3}f|H|^2-(\frac{\sqrt{30}}{6} f'-\frac{1}{6}f'')|\overline{\nabla} J_{\Sigma_t}|^2\leq 0.$$
We will find $f$ such that $\frac{\sqrt{30}}{6} f'-\frac{1}{6}f''\geq 0$. Noticing (\ref{e17}), it suffice to have
$$\frac{4}{3}f|H|^2-\frac{1}{12}(\sqrt{30}f'-f'')|H|^2\leq 0,$$ i.e,
$$\frac{4}{3}f-\frac{1}{12}(\sqrt{30}f'-f'')\leq 0.$$ Set $g=\frac{f'}{f}$, then $\frac{f''}{f}=g'+g^2$, then
it reduces to solve the
inequality
\begin{eqnarray*}
\frac{4}{3}-\frac{\sqrt{30}}{12}g+\frac{1}{12}g^2+\frac{1}{12}g' &\leq& 0,~~~~~~~~~~~~~~~~~~x\in [\frac{\sqrt{30}}{6}, 1],
\end{eqnarray*}
and
\begin{eqnarray*}
g'+g^2 &\geq& 0, ~~~~~~~~~~~~~~~~~~~~~~~~~~~~~~~~~~~~~~~~~~~~~~~~~~~~~x\in [\frac{\sqrt{30}}{6}, 1].
\end{eqnarray*}
 Assume that $g=-16x+b$, where $b$ is a constant, then it reduces to solve
$$4<-16x+b\leq \sqrt{30},~~~~~~~~~~~~~~~~~~x\in [\frac{\sqrt{30}}{6}, 1].$$ Therefore we choose $b$ such that
$$20\leq b\leq 8\frac{\sqrt{30}}{3}+\sqrt{30}.$$  We can choose $b=20$, i.e,
$$f=e^{-8x^2+20x}.$$ Thus
\begin{eqnarray*}
(\frac{\partial}{\partial t}-\Delta)\frac{|H|^2}{f(\cos\alpha)}
&\leq& 2\frac{f'\nabla\cos\alpha}{f}\cdot\nabla\frac{|H|^2}{f(\cos\alpha)}+\frac{9}{4}k\frac{|H|^2}{f(\cos\alpha)}.
\end{eqnarray*} This implies that
$$\frac{|H|^2}{f(\cos\alpha)}\leq e^{\frac{9k}{4}t}\frac{|H|^2}{f(\cos\alpha)}(0).$$ Since $\frac{\sqrt{30}}{6}\leq\cos\alpha\leq 1$,
$f(x)$ is bounded in $[\frac{\sqrt{30}}{6}, 1]$, we have
$$|H|^2\leq C_0 e^{\frac{9k}{4}t},$$ where $C_0$ depends only on $\max_{\Sigma_0}|H|^2$. Pinching inequality implies
 $|A|^2\leq  C_0 e^{\frac{9k}{4}t}+\frac{k}{2}.$ We finish the proof of the theorem. \hfill Q. E. D.

\begin{remark}
It was pointed out by Yang Liuqing that we could choose the linear function $f=x-\frac{3}{4}$ in the proof of the above theorem.
\end{remark}

\begin{theorem}\label{T2}
Under the assumption of Theorem \ref{L1}, the symplectic mean curvature flow converges to a holomorphic curve.
\end{theorem}

{\it Proof.} We can rewrite the evolution  equation of $\cos\alpha$
\begin{eqnarray*}
(\frac{\partial}{\partial t}-\Delta)\cos\alpha=|\overline{\nabla}J_{\Sigma_t}|^2\cos\alpha+\frac{3k}{2}\cos\alpha\sin^2\alpha.
\end{eqnarray*}  as
\begin{eqnarray}\label{sin}
(\frac{\partial}{\partial
t}-\Delta)\sin^2(\alpha/2)&=&-|\overline{\nabla}J_{\Sigma_t}|^2\cos\alpha-6k\sin^2(\alpha/2)\cos^2(\alpha/2)\cos\alpha
\\ &\leq& -c\sin^2(\alpha/2),
\end{eqnarray} where $c>0$ depends only on $k$ and the lower bound of  $\cos\alpha$. Applying the maximum principle, we get that
$\sin^2(\alpha/2)\leq e^{-ct}$. By Theorem \ref{T1} we know that the symplectic mean curvature flow exists for long time. Thus for
any $\varepsilon>0$, there exists $T$ such that as $t>T$, we have \begin{eqnarray}\label{e19}
\cos\alpha&\geq& 1-\varepsilon, \nonumber\\ \sin\alpha &\leq&2\varepsilon,\nonumber \\
|\nabla\cos\alpha|^2&\leq&2\varepsilon |\overline\nabla J_{\Sigma_t}|^2\leq 4\varepsilon |A|^2.
\end{eqnarray} Therefore,
\begin{eqnarray}\label{e20}
(\frac{\partial}{\partial t}-\Delta)\cos\alpha &\geq& \frac{1}{2}|H|^2\cos\alpha+\frac{3}{2}k\sin^2\alpha\cos\alpha\nonumber\\
&\geq&(\frac{3}{4}|A|^2-\frac{3k}{8})\cos\alpha +\frac{3}{2}k\sin^2\alpha\cos\alpha\nonumber\\ &\geq& \frac{3}{4}(1-\varepsilon)|A|^2-\frac{3k}{8}.
\end{eqnarray}

From (\ref{e8}) we see that
\begin{eqnarray*}
(\frac{\partial}{\partial t}-\Delta)|A|^2 \leq -2|\nabla A|^2+C_1|A|^4+C_2|A|^2+C_3,
\end{eqnarray*} where $C_1, C_2, C_3$ are constants that depend on the bounds of the curvature
tensor of $M$.

Let $p>1$ be a constant to be fixed later. For simplicity, we set $u=\cos\alpha$. Now we consider the function $\frac{|A|^2}{e^{pu}}$.
\begin{eqnarray*}
(\frac{\partial}{\partial t}-\Delta)\frac{|A|^2}{e^{pu}}&=& 2\nabla(\frac{|A|^2}{e^{pu}})\cdot\frac{\nabla e^{pu}}{e^{pu}}\\
&&+ \frac{1}{e^{2pu}}[e^{pu}(\frac{\partial}{\partial t}-\Delta)|A|^2-|A|^2(\frac{\partial}{\partial t}-\Delta)e^{pu}]\\ &\leq&
2p\nabla(\frac{|A|^2}{e^{pu}})\cdot\nabla u \\&&+\frac{1}{e^{2pu}}[e^{pu}(C_1|A|^4+C_2|A|^2+C_3)\\&&-p|A|^2e^{pu}[\frac{3}{4}(1-\varepsilon)|A|^2
-\frac{3k}{8}-p|\nabla u|^2]]
.
\end{eqnarray*}

Using (\ref{e19}) we obtain that,
\begin{eqnarray*}
(\frac{\partial}{\partial t}-\Delta)\frac{|A|^2}{e^{pu}} &\leq& 2p\nabla(\frac{|A|^2}{e^{pu}})\cdot\nabla u \\&&+
\frac{1}{e^{pu}}[(C_1-\frac{3p}{4}(1-\varepsilon)+4p^2\varepsilon)|A|^4+C_4|A|^2+C_3].
\end{eqnarray*}
Set $p^2=1/\varepsilon$, then
\begin{eqnarray*}
C_1-\frac{3p}{4}(1-\varepsilon)+4p^2\varepsilon &=&C_1-\frac{3}{4}\varepsilon^{-\frac{1}{2}}+\frac{3}{4}\varepsilon^{\frac{1}{2}}+4.
\end{eqnarray*} As $t$ is sufficiently large, i.e. $\varepsilon$ is sufficiently close to $0$, we have
\begin{eqnarray*}
(C_1-\frac{3}{4}\varepsilon^{-\frac{1}{2}}+\frac{3}{4}\varepsilon^{\frac{1}{2}}+4)\leq -1.
\end{eqnarray*} So,
\begin{eqnarray*}
(\frac{\partial}{\partial t}-\Delta)\frac{|A|^2}{e^{pu}} &\leq& 2p\nabla(\frac{|A|^2}{e^{pu}})\cdot\nabla u
-\frac{|A|^4}{e^{pu}}+C_4\frac{|A|^2}{e^{pu}}+\frac{C_3}{e^{pu}}\\ &\leq& 2p\nabla(\frac{|A|^2}{e^{pu}})\cdot\nabla u
-\frac{|A|^4}{e^{2pu}}+C_4\frac{|A|^2}{e^{pu}}+\frac{C_3}{e^{pu}}
\end{eqnarray*}

Applying the maximum principle for parabolic equations, we conclude that $\frac{|A|^2}{e^{pu}}$ is uniformly bounded, thus
$|A|^2$ is also uniformly bounded. Thus $F(\cdot, t)$ converges to $F_\infty$ in $C^2$ as $t\to\infty$. Since
$\sin^2(\alpha/2)\leq e^{-ct},$ we have $\cos\alpha\equiv 1$ at infinity. Thus the limiting surface $F_\infty$ is
a holomorphic curve.

\hfill Q. E. D.

\section{Another pinching estimate}

In this section, we will derive another pinching condition.

\begin{theorem}\label{T3}
Suppose $M$ is a K\"ahler surface with constant holomorphic sectional curvature $k>0$ and $\Sigma$ is a symplectic surface in $M$. Assume that
$|A|^2\leq \frac{2}{3}|H|^2+\frac{4}{5}k\cos\alpha$ and $\cos\alpha\geq\frac{251}{265}$ holds on the initial surface, then it remains true along the symplectic mean curvature flow. Furthermore, the symplecitic mean curvature flow exists for long time and converges to a holomorphic curve at infinity.
\end{theorem}

{\it Proof.} From (\ref{e33}), (\ref{e29}), we see that,
\begin{eqnarray*}
\frac{\partial}{\partial t}|A|^{2}&\leq&
     \Delta |A|^{2} -2|\nabla A|^{2}-k|A|^2-\frac{k}{2}(3\cos^2\alpha+1)|A|^2+2k|H|^2\\ &&
+2\sum_{\alpha,\beta,i,j}(\sum_{k}(h^{\alpha}_{ik}h^{\beta}_{jk}-h^{\alpha}_{jk}h^{\beta}_{ik}))^{2}
       +2\sum_{\alpha,\beta}(\sum_{ij}h^{\alpha}_{ij}h^{\beta}_{ij})^{2},
\end{eqnarray*}
and
\begin{eqnarray*}
\frac{\partial}{\partial t}|H|^{2}&=&
    \Delta |H|^{2}-2|\nabla H|^{2}+3k|H|^2-\frac{k}{2}(3\cos^2\alpha+1)|H|^2\\&&
               +2\sum_{ij}(\sum_{\alpha} H^{\alpha}h^{\alpha}_{ij})^2.
\end{eqnarray*}
Set $Q=|A|^2-\frac{2}{3}|H|^2-bk\cos\alpha$. Then we have,\allowdisplaybreaks
\begin{eqnarray}\label{e30}
\frac{\partial}{\partial t}Q &\leq&\Delta Q-2(|\nabla A|^{2}-\frac{2}{3}|\nabla H|^{2})-k|A|^2\nonumber\\ &&
-\frac{k}{2}(3\cos^2\alpha+1)(|A|^2-\frac{2}{3}|H|^2)\nonumber\\ &&-bk
(|\overline\nabla_{\Sigma_t}J|^2\cos\alpha+\frac{3}{2}k\cos\alpha\sin^2\alpha)\nonumber\\&&
+2\sum_{\alpha,\beta,i,j}(\sum_{k}(h^{\alpha}_{ik}h^{\beta}_{jk}-h^{\alpha}_{jk}h^{\beta}_{ik}))^{2}
       +2\sum_{\alpha,\beta}(\sum_{ij}h^{\alpha}_{ij}h^{\beta}_{ij})^{2}\nonumber\\&&-\frac{4}{3}
\sum_{ij}(\sum_{\alpha} H^{\alpha}h^{\alpha}_{ij})^2\nonumber\\ &\leq&  \Delta Q-2(|\nabla A|^{2}-\frac{2}{3}|\nabla H|^{2})-\frac{k}{2}(3\cos^2\alpha+1)Q
\nonumber\\ &&-2bk^2\cos\alpha-k|A|^2-\frac{bk}{2}\cos\alpha |H|^2
\nonumber\\&&
+2\sum_{\alpha,\beta,i,j}(\sum_{k}(h^{\alpha}_{ik}h^{\beta}_{jk}-h^{\alpha}_{jk}h^{\beta}_{ik}))^{2}
       +2\sum_{\alpha,\beta}(\sum_{ij}h^{\alpha}_{ij}h^{\beta}_{ij})^{2}\nonumber\\&&-\frac{4}{3}
\sum_{ij}(\sum_{\alpha} H^{\alpha}h^{\alpha}_{ij})^2.\nonumber
\end{eqnarray}
By an argument similar to the one used in the proof of Theorem \ref{L1}, at the point $H\neq 0$ we can get that
\begin{eqnarray*}
\frac{\partial}{\partial t}Q &\leq&\Delta Q+\frac{9k^2}{2}\cos^2\alpha\sin^2\alpha
-\frac{k}{2}(3\cos^2\alpha+1)Q\\&&-2bk^2\cos\alpha-k|A|^2-\frac{bk}{2}\cos\alpha |H|^2\\ &&+
2 |\mathring{h}_3|^4+2|\mathring{h}_4|^4+\frac{2}{3}|\mathring{h}_3|^2|H|^2-\frac{1}{6}|H|^4+8|\mathring{h}_3|^2|\mathring{h}_4|^2.
\end{eqnarray*} Sine $|H|^2=6(|\mathring{h}_3|^2+|\mathring{h}_4|^2-Q-bk\cos\alpha)$, putting it into the above inequality we obtain that
\begin{eqnarray}\label{e31}
\frac{\partial}{\partial t}Q &\leq&\Delta Q-6Q^2\nonumber
\\&&+[8|\mathring{h}_3|^2+12|\mathring{h}_4|^2-9bk\cos\alpha+3k-\frac{k}{2}(3\cos^2\alpha+1)]Q\nonumber\\&&
+\frac{9k^2}{2}\cos^2\alpha\sin^2\alpha+bk^2\cos\alpha-3b^2k^2\cos^2\alpha\nonumber\\&&-4|\mathring{h}_4|^2+k(5b\cos\alpha-4)|\mathring{h}_3|^2
+k(9b\cos\alpha-4)|\mathring{h}_4|^2\nonumber\\ &\leq&
\Delta Q-6Q^2\nonumber
\\&&+[8|\mathring{h}_3|^2+12|\mathring{h}_4|^2-9bk\cos\alpha+3k-\frac{k}{2}(3\cos^2\alpha+1)]Q\nonumber\\&&
+k(5b\cos\alpha-4)|\mathring{h}_3|^2-(2|\mathring{h}_4|^2-\frac{k}{4}(9b\cos\alpha-4))^2+\frac{k^2}{16}(9b\cos\alpha-4)^2\nonumber\\
&&+\frac{9k^2}{2}\cos^2\alpha\sin^2\alpha+bk^2\cos\alpha-3b^2k^2\cos^2\alpha\nonumber\\ &\leq&
\Delta Q-6Q^2\nonumber
\\&&+[8|\mathring{h}_3|^2+12|\mathring{h}_4|^2-9bk\cos\alpha+3k-\frac{k}{2}(3\cos^2\alpha+1)]Q\nonumber\\&&
+k(5b\cos\alpha-4)|\mathring{h}_3|^2-(2|\mathring{h}_4|^2-\frac{k}{4}(9b\cos\alpha-4))^2\nonumber\\&&+\frac{k^2}{16}
(33b^2\cos^2\alpha-56b\cos\alpha+16+72\cos^2\alpha\sin^2\alpha)
\end{eqnarray} Now we need choose $b$ and  the lower bound of $\cos\alpha$ such that  $5b\cos\alpha-4\leq 0$ and
$$33b^2\cos^2\alpha-56b\cos\alpha+16+72\cos^2\alpha\sin^2\alpha\leq 0.$$ First we choose $b=\frac{4}{5}$, then we need
\begin{eqnarray}\label{e40}
33\times\frac{16}{25}\cos^2\alpha-56\times\frac{4}{5}\cos\alpha+16+72\cos^2\alpha\sin^2\alpha\leq 0.
\end{eqnarray}
Assume that $\cos\alpha\geq\delta$. If
\begin{eqnarray}\label{e41}
33\times\frac{16}{25}\cos^2\alpha-56\times\frac{4}{5}\cos\alpha+16+72(1-\cos^2\alpha)\leq 0
\end{eqnarray} holds, then (\ref{e40}) holds. Solving (\ref{e41}), we get that
\begin{equation}\label{e42}
\delta\geq\frac{251}{265}.
\end{equation}

At the point $H=0$, using (\ref{e35}), we have
\begin{eqnarray*}
(\frac{\partial}{\partial t}-\Delta)Q &\leq&\frac{9}{2}k^2\sin^2\alpha\cos^2\alpha-\frac{k}{2}(3\cos^2\alpha+1)Q\\
&&-2bk^2\cos\alpha-k|A|^2+3|A|^4.
\end{eqnarray*}
Putting $|A|^2=Q+bk\cos\alpha$ into the above inequality, we get that,
\begin{eqnarray*}
(\frac{\partial}{\partial t}-\Delta)Q &\leq& [3(|A|^2-bk\cos\alpha)+6bk\cos\alpha-k-\frac{k}{2}(3\cos^2\alpha+1)]Q \\ &&
+\frac{9}{2}k^2\sin^2\alpha\cos^2\alpha-3bk^2\cos\alpha+3b^2k^2\cos^2\alpha.
\end{eqnarray*}
Choose $b=\frac{4}{5}$ and assume that $\cos\alpha\geq\delta$, then we need
\begin{eqnarray*}
\frac{9}{2}(1-\cos\alpha^2)-\frac{12}{5}\cos\alpha+\frac{48}{25}\cos^2\alpha\leq 0.
\end{eqnarray*}
Solving it we get that
\begin{eqnarray}\label{e43}
\delta\geq\frac{121}{129}.
\end{eqnarray}

Compare (\ref{e42}) and (\ref{e43}), we choose $\delta\geq\frac{251}{265}$ and $b=\frac{4}{5}$.

The global existence and convergence of the symplectic mean curvature flow can be proved in a similar manner as the one used in the proof
of Theorem \ref{T1} and Theorem \ref{T2}.
\hfill Q. E. D.

\begin{remark}
It is clear that the pinching condition $|A|^2\leq \frac{2}{3}|H|^2+\frac{4}{5}k\cos\alpha$
is better than the condition that $|A|^2\leq\frac{2}{3}|H|^2+\frac{1}{2}k$. On the other hand, the condition $\cos\alpha\geq\frac{\sqrt{30}}{6}$
is better than $\cos\alpha\geq\frac{251}{265}$.
\end{remark}

\end{document}